\documentclass[12pt,reqno,twoside]{article}%draft
\usepackage{amsmath,amsthm,amsfonts,amssymb,setspace,verbatim,cite,fancyhdr}

\usepackage[english]{babel}
\usepackage[utf8]{inputenc}
\usepackage{times}

\usepackage[T1]{fontenc}
\usepackage{mathrsfs}

\let\oldnsubseteq\nsubseteq 
\usepackage{mathabx} 
\let\nsubseteq\oldnsubseteq

\usepackage{enumitem}
\usepackage{accents}
\usepackage{fge}

\usepackage[top=15mm,bottom=13mm,left=22.8mm,right=22.8mm,includeheadfoot,headsep=8mm,footnotesep=5mm]{geometry} 

\usepackage[hypertexnames=false]{hyperref} 
\hypersetup{
pdfencoding=auto,
pageanchor=true,
plainpages=true,
pdftoolbar=true,
pdfmenubar=true,
pdfpagemode=UseOutlines,%(show bookmarks);UseNone;UseThumbs (show thumbnails);FullScreen;UseOC (PDF % 1.5);UseAttachments(PDF 1.6)
pdffitwindow=false,
pdfstartview={FitH},
pdfnewwindow=true,
colorlinks=true,
linkcolor=black,
citecolor=black,
filecolor=black,
urlcolor=black,
linktoc=all,
breaklinks=true,
}

\usepackage{mdframed,xcolor}
\usepackage{xpatch}

\makeatletter
\xpatchcmd{\endmdframed}
  {\aftergroup\endmdf@trivlist\color@endgroup}
  {\endmdf@trivlist\color@endgroup\@doendpe}
  {}{}
\makeatother

\mdfdefinestyle{defStyle}
{linewidth=0.4mm,
skipabove=0.5\baselineskip,
skipbelow=0.9\baselineskip,
splittopskip=\topskip,
innertopmargin=0pt,
theoremspace={},
theoremtitlefont=\normalfont,
}

\newmdtheoremenv[linewidth=0.4mm,skipabove=0.5\baselineskip,skipbelow=0.5\baselineskip,splittopskip=\topskip,%
innertopmargin=0pt]{theorem}{Theorem}[section]
\newmdtheoremenv[linewidth=0.4mm,skipabove=0.5\baselineskip,skipbelow=0.5\baselineskip,splittopskip=\topskip,%
innertopmargin=0pt]{lemma}[theorem]{Lemma}
\newmdtheoremenv[linewidth=0.4mm,skipabove=0.5\baselineskip,skipbelow=0.5\baselineskip,splittopskip=\topskip,%
innertopmargin=0pt]{corollary}[theorem]{Corollary}
\newmdtheoremenv[linewidth=0.4mm,skipabove=0.5\baselineskip,skipbelow=0.5\baselineskip,splittopskip=\topskip,%
innertopmargin=0pt]{proposition}[theorem]{Proposition}
\newtheoremstyle{myexercise}{-\topsep}{}{\normalfont}{}{\bfseries}{}{.5em}{}
\theoremstyle{myexercise}

\theoremstyle{definition}
\newmdtheoremenv[linewidth=0.4mm,skipabove=0.5\baselineskip,skipbelow=0.5\baselineskip,splittopskip=\topskip,%
innertopmargin=0pt]{hypothesis}[theorem]{Hypothesis}
\newmdtheoremenv[linewidth=0.4mm,skipabove=0.5\baselineskip,skipbelow=0.5\baselineskip,splittopskip=\topskip,%
innertopmargin=0pt]{definition}[theorem]{Definition}

\usepackage{thmtools}

\declaretheorem[style=definition,name=Remark,sibling=theorem]{remark}

\numberwithin{equation}{section}

\makeatletter
\def\th@plain{%
  \thm@notefont{\bfseries}% same as heading font
  \itshape % body font
}

\def\th@definition{%
  \normalfont % body font
  \thm@notefont{\bfseries}% same as heading font
}
\makeatother

\renewcommand{\thefootnote}{\fnsymbol{footnote}}

\newcommand{\makepapertitle}{
\pdfbookmark[1]{Title page}{Title_page}
\thispagestyle{empty}\vspace*{8mm}
\begin{spacing}{1.1}
\LARGE\MakeUppercase{\papertitle}
\end{spacing}\vspace*{10mm}
}

\newcommand{\paperfirstauthor}{%
{\large\sc\firstauthor}\\[2mm] %\FootCorAuth%\footnote{\firstthanks}} \\[2mm]
{\rm\firstaddress \\ 
\firstemail}\\[8mm]
}

\newcommand{\makeabstract}{%
\begin{minipage}{0.9\textwidth}
{\small {\sc Abstract.}
 \paperabstract
}
\end{minipage}\vfill
}

\newcommand{\MakeFirstPageOneAuthor}{
\begin{center}
  \makepapertitle
  \paperfirstauthor
  \makeabstract
  \begin{minipage}[t]{0.3\textwidth}
   \raggedleft {\bf Date (final version):} 
  \end{minipage}\hspace{0.01\textwidth}
  \begin{minipage}[t]{0.6\textwidth}
    \today\\ 
    (submitted on June 04, 2024;\\ accepted on December 20, 2024)
  \end{minipage}\\[4mm]
\noindent%
  \begin{minipage}[t]{0.3\textwidth}
    \raggedleft {\bf Running head:} 
  \end{minipage}\hspace{0.01\textwidth}
  \begin{minipage}[t]{0.6\textwidth}
    \runninghead
  \end{minipage}\\[4mm]
\noindent%
  \begin{minipage}[t]{0.3\textwidth}
    \raggedleft {\bf Submitted to:} 
  \end{minipage}\hspace{0.01\textwidth}
  \begin{minipage}[t]{0.6\textwidth}
    \submittedto
  \end{minipage}\\[4mm]
\bigskip
\vfill
\end{center}

\renewcommand{\thefootnote}{}
\footnotetext[2]{2020 {\it Mathematics Subject Classification:} \thesubjclass}
\footnotetext[3]{{\it Key words and phrases:} \thekeywords}
\setcounter{footnote}{0}
\renewcommand{\thefootnote}{{\bf\,\alph{footnote}\alph{footnote}\alph{footnote}\,}}

\setcounter{page}{0}
\newpage
}

\makeatletter
\newcommand{\rightorleftmark}{%
  \begingroup\protected@edef\x{\rightmark}%
  \ifx\x\@empty
    \endgroup\leftmark
  \else
    \endgroup\rightmark
  \fi}
\makeatother

\newcommand{\papertitle}%
{The Friedrichs extension of a class of discrete symplectic systems}

\newcommand{\runninghead}%
{Friedrichs extension and discrete symplectic systems}

\newcommand{\firstauthor}%
{Petr Zem{\'{a}}nek}
\newcommand{\firstauthorhead}%
{P. Zem{\'{a}}nek}

\newcommand{\firstaddress}
{Department of Mathematics and Statistics, Faculty of Science, Masaryk University \\
Kotl{\'{a}}{\v{r}}sk{\'{a}} 2, CZ-61137 Brno, Czech Republic}

\newcommand{\firstemail}%
{E-mail: zemanekp@math.muni.cz}

\newcommand{\paperabstract}%
{The Friedrichs extension of minimal linear relation being bounded below and associated with the discrete symplectic
system with a special linear dependence on the spectral parameter is characterized by using recessive solutions. This
generalizes a similar result obtained by Došlý and Hasil for linear operators defined by infinite banded matrices
corresponding to even-order Sturm--Liouville difference equations and, in a certain sense, also results of Marletta and
Zettl or Šimon Hilscher and Zemánek for singular differential operators.}

\newcommand{\thekeywords}%
{Discrete symplectic system; Friedrichs extension; minimal linear relation; recessive solution.}

\newcommand{\thesubjclass}%
{{\it Primary\/} 47A06; {\it Secondary\/} 47A20; 47B39; 39A06; 39A12.}

\newcommand{\submittedto}%
{Journal of Spectral Theory}

\DeclareMathAccent{\wwtilde}{\mathord}{largesymbols}{"65}
\DeclareMathSymbol{\widetildesym}{\mathord}{largesymbols}{"65}
\newcommand\lowerwidetildesym{%
  \text{\smash{\raisebox{-1.3ex}{%
    $\widetildesym$}}}}
\newcommand\wtilde[1]{%
  \mathchoice
    {\accentset{\displaystyle\lowerwidetildesym}{#1}}
    {\accentset{\textstyle\lowerwidetildesym}{#1}}
    {\accentset{\scriptstyle\lowerwidetildesym}{#1}}
    {\accentset{\scriptscriptstyle\lowerwidetildesym}{#1}}
}

\newcommand\lowerwidetildesymW{%
  \text{\smash{\raisebox{-1.35ex}{$\widetildesym$}}}}
\newcommand\wtildeW[1]{%
    {\accentset{\scalebox{1.2}{\lowerwidetildesymW}}{#1}}
}

\DeclareMathSymbol{\widehatsym}{\mathord}{largesymbols}{"62}
\newcommand\lowerwidehatsym{%
  \text{\smash{\raisebox{-1.3ex}{%
    $\widehatsym$}}}}
\newcommand\what[1]{%
  \mathchoice
    {\accentset{\displaystyle\lowerwidehatsym}{#1}}
    {\accentset{\textstyle\lowerwidehatsym}{#1}}
    {\accentset{\scriptstyle\lowerwidehatsym}{#1}}
    {\accentset{\scriptscriptstyle\lowerwidehatsym}{#1}}
}

\DeclareMathAlphabet{\mthdtcl}{U}{dutchcal}{m}{n}
\DeclareMathAlphabet{\mathpzc}{OT1}{pzc}{m}{it}
\DeclareMathAlphabet{\msfsl}{U}{eus}{m}{n}

\newcommand{\e}{\mathrm{e}}

\newcommand{\Ac}{\mathcal{A}}

\newcommand{\Bc}{\mathcal{B}}

\newcommand{\Cc}{\mathcal{C}}

\newcommand{\Dc}{\mathcal{D}}

\newcommand{\Fc}{\mathcal{F}}

\newcommand{\Ic}{\mathcal{I}}

\newcommand{\Jc}{\mathcal{J}}

\newcommand{\Sc}{\mathcal{S}}

\newcommand{\Vc}{\mathcal{V}}

\newcommand{\Wc}{\mathcal{W}}
\newcommand{\tWc}{\wtildeW{\Wc}}

\newcommand{\mL}{\mathscr{L}}

\newcommand{\Cbb}{\mathbb{C}}

\newcommand{\Nbb}{\mathbb{N}}
\newcommand{\Rbb}{\mathbb{R}}
\newcommand{\Sbb}{\mathbb{S}}

\newcommand{\Zbb}{\mathbb{Z}}

\newcommand{\la}{\lambda}

\newcommand{\La}{\Lambda}
\newcommand{\Ups}{\Upsilon}

\newcommand{\bla}{\bar{\la}}

\newcommand{\Ps}{\Psi}

\newcommand{\Th}{\Theta}

\newcommand{\si}{\sigma}

\newcommand{\stm}{\hspace*{0.2mm}\fgebackslash\hspace*{0.3mm}}

\usepackage{scalerel}
\usepackage{stackengine}

\stackMath

\newcommand{\dbtilde}[1]{\accentset{\approx}{#1}}
                       
\newlength{\dhatheight}
\newcommand{\dbhat}[1]{%
    \settoheight{\dhatheight}{\ensuremath{\hat{#1}}}%
    \addtolength{\dhatheight}{-0.35ex}%
    \what{\vphantom{\rule{1pt}{\dhatheight}}%
    \smash{\what{#1}}}}                       

\newcommand{\dbhatsmall}[1]{%
    \settoheight{\dhatheight}{\ensuremath{\hat{#1}}}%
    \addtolength{\dhatheight}{-0.35ex}%
    \what{\vphantom{\rule{1pt}{0.7\dhatheight}}%
    \smash{\what{#1}}}}

\newcommand{\pzcH}{\mathpzc{H}}

\newcommand{\pzcT}{\mathpzc{T}}

\newcommand{\pzcX}{\mathpzc{X}}
\newcommand{\pzcY}{\mathpzc{Y}}

\newcommand{\tf}{\tilde{f}}
\newcommand{\tfjm}{\tilde{f}^{[j,m]}}

\newcommand{\fclass}{[f]}
\newcommand{\pzcf}{\mathpzc{f}}

\newcommand{\ttf}{\dbtilde{f}}
\newcommand{\ttfjm}{\ttf^{[j,m]}}

\newcommand{\pzcg}{\mathpzc{g}}

\newcommand{\tu}{\tilde{u}}

\newcommand{\tU}{\wtilde{U}}
\newcommand{\hU}{\what{U}}

\newcommand{\hhU}{\dbhat{U}}
\newcommand{\hhhU}{\dbhatsmall{U}}
\newcommand{\Um}{U^{[m]}}
\newcommand{\ujm}{u^{[j,m]}}
\newcommand{\tujm}{\tu^{[j,m]}}

\newcommand{\wclass}{[w]}

\newcommand{\tx}{\tilde{x}}

\newcommand{\tX}{\wtilde{X}}
\newcommand{\hX}{\what{X}}

\newcommand{\hhX}{\dbhat{X}}
\newcommand{\Xm}{X^{[m]}}
\newcommand{\xjm}{x^{[j,m]}}
\newcommand{\txjm}{\tx^{[j,m]}}

\newcommand{\pzcy}{\mathpzc{y}}

\newcommand{\tz}{\tilde{z}}
\newcommand{\zclass}{[z]}
\newcommand{\tZ}{\wtilde{Z}}

\newcommand{\hhZ}{\dbhat{Z}}

\newcommand{\ttz}{\dbtilde{z}}
\newcommand{\ttzjm}{\ttz^{[j,m]}}
\newcommand{\ttzj}{\ttz^{[j]}}
\newcommand{\hz}{\hat{z}}
\newcommand{\hZ}{\what{Z}}
\newcommand{\pzcz}{\mathpzc{z}}

\newcommand{\Zm}{Z^{[m]}}
\newcommand{\zjm}{z^{[j,m]}}
\newcommand{\tzjm}{\tz^{[j,m]}}
\newcommand{\tzj}{\tz^{[j]}}

\newcommand{\ltp}{\ell^{\hspace{0.2mm}2}_{\Ps}}

\newcommand{\tltp}{\tilde{\ell}^{\hspace{0.3mm}2}_{\Ps}}

\newcommand{\Tmax}{T_{\mathrm{max}}}
\newcommand{\Tmin}{T_{\mathrm{min}}}

\newcommand{\sZbb}{{\scriptscriptstyle{\Zbb}}}
\newcommand{\IzD}{\Ic_\sZbb^{\scriptscriptstyle{\rm D}}}
\newcommand{\oinftyZ}{[0,\infty)_\sZbb}

\newcommand{\mmatrix}[1]{\begin{pmatrix} #1
  \end{pmatrix}}
\newcommand{\msmatrix}[1]{\left(\begin{smallmatrix} #1
  \end{smallmatrix}\right)}  

\newcommand{\qtextq}[1]{\quad\text{#1}\quad}
\newcommand{\qtext}[1]{\quad\text{#1 }\ }

\DeclareMathOperator{\re}{Re}

\DeclareMathOperator{\rank}{rank}
\DeclareMathOperator{\dom}{dom}

\usepackage{mathtools,xparse}
\renewcommand{\.}{\hspace*{0.1 em}}
\DeclarePairedDelimiter\xnorm{\lVert}{\rVert}
\NewDocumentCommand{\norm}{som}
 {\IfBooleanTF{#1}
   {\xnorm*{#3}}
   {\IfNoValueTF{#2}
     {\mathopen{|\mkern-.8mu|}\.#3\.\mathclose{|\mkern-.8mu|}}
     {\xnorm[#2]{\.#3\.}}%
   }
 }
\DeclarePairedDelimiter\xinner{\langle}{\rangle}

\NewDocumentCommand{\xinnr}{som}
 {\IfBooleanTF{#1}
   {\xinner*{#3}}
   {\IfNoValueTF{#2}
     {\mathopen{\langle}\.#3\.\mathclose{\rangle}}
     {\xinner[#2]{\.#3\.}}%
   }
 }

\makeatletter
\def\inner{\@ifnextchar[{\@INNwith}{\@INNwithout}}
\def\@INNwith[#1]#2#3{\xinnr[#1]{#2,#3}}
\def\@INNwithout#1#2{\xinnr{#1,#2}}

\def\innerP{\@ifnextchar[{\@INNPwith}{\@INNPwithout}}
\def\@INNPwith[#1]#2#3{\xinnr[#1]{#2,#3}_\Ps}
\def\@INNPwithout#1#2{\xinnr{#1,#2}_\Ps}

\def\innerPN{\@ifnextchar[{\@INNPNwith}{\@INNPNwithout}}
\def\@INNPNwith[#1]#2#3{\xinnr[#1]{#2,#3}_{\Ps,N}}
\def\@INNPNwithout#1#2{\xinnr{#1,#2}_{\Ps,N}}

\def\abs{\@ifnextchar[{\@awith}{\@awithout}}
\def\@awith[#1]#2{{#1|}#2\.{#1|}}
\def\@awithout#1{|#1\.|}
 
\def\normS{\@ifnextchar[{\@Nwith}{\@Nwithout}}
\def\@Nwith[#1]#2{\norm[#1]{#2}_\si}
\def\@Nwithout#1{\norm{#1}_\si}

\def\normE{\@ifnextchar[{\@NEwith}{\@NEwithout}}
\def\@NEwith[#1]#2{\norm[#1]{#2}_2}
\def\@NEwithout#1{\norm{#1}_2}

\def\normA{\@ifnextchar[{\@NAwith}{\@NAwithout}}
\def\@NAwith[#1]#2{\norm[#1]{#2}_1}
\def\@NAwithout#1{\norm{#1}_1}

\def\normP{\@ifnextchar[{\@NPwith}{\@NPwithout}}
\def\@NPwith[#1]#2{\norm[#1]{#2}_{\Ps}}
\def\@NPwithout#1{\norm{#1}_{\Ps}}

\def\normW{\@ifnextchar[{\@NWwith}{\@NWwithout}}
\def\@NWwith[#1]#2{\norm[#1]{#2}_{\Wc}}
\def\@NWwithout#1{\norm{#1}_{\Wc}}

\def\normtW{\@ifnextchar[{\@NtWwith}{\@NtWwithout}}
\def\@NtWwith[#1]#2{\norm[#1]{#2}_{\scalebox{0.7}{$\tWc$}}}
\def\@NtWwithout#1{\norm{#1}_{\scalebox{0.6}{$\tWc$}}}

\newcommand{\Sla}[1]{\text{\rm(S$_{#1}$})}
\newcommand{\Slaf}[2]{\text{\rm(S$_{#1}^{#2}$)}}

\makeatletter
\newcommand{\ltxlabel}{\ltx@label}
\makeatother
\newcounter{GatherItemCounter}

\usepackage{ifdraft}
\ifdraft{\overfullrule26pt}{\overfullrule0pt}
\ifdraft{\usepackage[nomsgs]{refcheck}}{}
\ifdraft{\geometry{showframe,showcrop}}{}

\hypersetup{
pdfencoding=unicode,
pdftitle={\papertitle},
pdfauthor={\firstauthor},
pdfsubject={Department of Mathematics and Statistics, Faculty of Science, Masaryk University},
pdfkeywords={\thekeywords},
pdfsubject=\paperabstract,
pdfencoding=auto
}

\begin{document}

%%%%%%%%%%%%%%%%%%%%%%%%%%%%%%%%%%%%%%%%%%%%%% FIRST PAGE %%%%%%%%%%%%%%%%%%%%%%%%%%%%%%%%%%%%%%%%%%%%%%%%%%%%%%%%%%%%

\MakeFirstPageOneAuthor

%%%%%%%%%%%%%%%%%%%%%%%%%%%%%%%%%%%%%%%%%%%% SECTION %%%%%%%%%%%%%%%%%%%%%%%%%%%%%%%%%%%%%%%%%%%%%%%%%%%%%%%%%%%%%%%%%

\section{Introduction}\label{S:intro}

Qualitative properties of operators or (more generally) linear relations can be investigated in various ways, including 
a structure of their spectrum, boundary triplets, or a description of their self-adjoint extensions with a focus on some 
particular cases. Especially the Friedrichs extension belongs to the very traditional topics, and it has attracted more 
attention again in recent years, see e.g. 
\cite{qB.gW.aZ23,cY.hS21,sZ.hS.cY23,zZ.qK18,zZ.jQ.jS22,zZ.jQ.jS23,sY.jS.aZ15,qB.gW.aZ22:OaM,qB.gW.aZ22:JMAA,ogS18,kS23,
kS21,zS.zT21}. Therefore, in the present paper, we aim to characterize the (domain of the) Friedrichs 
extension of the minimal linear relation determined by the linear mapping 
 \begin{equation*}%\label{E:mL.def}
  \mL(z)_k\coloneq\Jc(z_k-\Sc_k\,z_{k+1})
 \end{equation*}
acting on a weighted space of $2n$-vector valued square summable sequences $\ltp$ with respect to the weight matrices 
$\Ps_k$ on the unbounded discrete interval $\oinftyZ\coloneq [0,\infty)\cap\Zbb$, where the coefficients are 
$2n\times2n$ complex-valued matrices satisfying
 \begin{equation}\label{E:Sc.Ps.intro}
  \Sc_k^*\Jc\Sc_k=\Jc
  \qtextq{and}
  \Ps_k^*\,\Jc\,\Ps_k=\Ps_k\,\Jc\,\Ps_k=0 \qtext{for all $k\in\oinftyZ$}
 \end{equation}
with the superscript $*$ denoting the conjugate transpose and $\Jc$ standing for the $2n\times2n$ orthogonal and 
skew-symmetric matrix  
 \begin{equation}\label{E:matrix.J.def}
  \Jc=\Jc_{2n}\coloneq\mmatrix{0 & I_n\\ -I_n & 0}.
 \end{equation}
The first equality in~\eqref{E:Sc.Ps.intro} means that $\Sc_k$ is symplectic for all $k$ and the mapping $\mL$ is 
closely related to the (nonhomogeneous) {\it time-reversed discrete symplectic system} on the half-line, because the 
relation $\mL(z)=\la\.\Ps\.z+\Ps\.f$ with arbitrary $\la\in\Cbb$ is equivalent to
 \begin{equation*}\label{E:Sla}\tag{S$_\la^f$}
  z_{k}(\la)=(\Sc_k+\la\.\Vc_k)\,z_{k+1}(\la)-\Jc\.\Ps_k\.f_k, \quad k\in\oinftyZ,
 \end{equation*}
where $\Vc_k=-\Jc\,\Ps_k\,\Sc_k$ is such that 
 \begin{equation*}\label{E:Sc.Vc.asmpt}
  \Vc_k^*\Jc\Sc_k\ \text{ is Hermitian,} \qtextq{and} 
  \Vc_k^*\Jc\,\Vc_k=0 \qtext{for all $k\in\oinftyZ$.}
 \end{equation*}
The associated {\it minimal linear relation} can be written as
 \begin{equation}\label{E:Tmin.explicit}
  \Tmin=\big\{\{\zclass,\fclass\}\in\Tmax\mid z_0=0=\lim_{k\to\infty}z_k^*\.\Jc\.w_k\ \text{ for all }
              \wclass\in\dom\Tmax\big\},
 \end{equation}
which is a restriction of the {\it maximal linear relation} given by
 \begin{equation}\label{E:Tmax.def}
  \Tmax\coloneq \big\{\{\zclass,\fclass\}\mid \text{there exists $u\in\zclass$ such that }
                                                     \mL(u)_k=\Ps_k\,f_k\ \text{ for all } k\in\oinftyZ\big\},
 \end{equation}
where $\zclass,\fclass$ stand for equivalence classes in $\ltp$. Actually, the proper definition of the minimal linear 
relation guarantees that its adjoint relation is $\Tmax$, i.e., it holds $\Tmin^*=\Tmax$ as shown 
in~\cite[Theorem~5.10]{slC.pZ15}. These relations and square summability of solutions of discrete symplectic systems 
were thoroughly studied by the author and his collaborators 
in~\cite{slC.pZ15,pZ.slC16,pZ.slC22,rSH.pZ13,rSH.pZ14:ICDEA,rSH.pZ14:JDEA,rSH.pZ15:MN,slC.pZ10}. Now, we turn our
attention to the Friedrichs extension $T_F$ of $\Tmin$, which is defined as a self-adjoint extension of $\Tmin$ being 
bounded below by the same lower bound as $\Tmin$. However, even though we speak of an extension, this linear relation 
$T_F$ will be expressed as a~{\it restriction} of $\Tmax$, which consists of pairs satisfying a~zero boundary condition 
at $k=0$ and a~specific limit condition at $\infty$ determined by recessive solutions of~\Slaf{\la}{f} with $f\equiv0$, 
i.e., of the (homogeneous) system
  \begin{equation*}\label{E:Sla.noref}\tag{S$_\la$}
  z_{k}(\la)=(\Sc_k+\la\.\Vc_k)\,z_{k+1}(\la), \quad k\in\oinftyZ,
 \end{equation*}
see Theorem~\ref{T:main} for a~precise formulation. Our main result relies on several facts from the theory of discrete 
symplectic systems. The first is a connection between the boundedness from below of $\Tmin$ and the existence of 
recessive solutions of~\eqref{E:Sla.noref}. The second crucial ingredient is the recessive solution of 
\eqref{E:Sla.noref} {\it per se}, because its properties imply the square integrability and, roughly speaking, presence 
in the domain of $T_F$. Here we should emphasize that recessive solutions are defined through the behavior of their 
first $n$ components of the $2n$-vector-valued solutions, which naturally leads to the restriction that we consider only 
the case when the weight matrices $\Ps_k$ have for all $k\in\oinftyZ$ the very special block structure
 \begin{equation}\label{E:Ps.block.def}
  \Ps_k\coloneq\mmatrix{\Wc_k & 0\\ 0 & 0}
 \end{equation}
with $\Wc_k=\Wc_k^*>0$ being $n\times n$ matrices. Finally, we utilize the characterization of all self-adjoint 
extensions of $\Tmin$ established in~\cite[Theorem~3.3 and~Remark~3.4]{pZ.slC16} and give precisely $d$ boundary 
conditions determining the Friedrichs extension, see also Theorem~\ref{T:intro.sadj.GKN.Sla} and 
Remark~\ref{R:overdetermined}\ref{R:overdetermined.i}.

The origin of the study of the concept nowadays known as the Friedrichs extension can be traced back to von 
Neumann. He showed that for any Hermitian linear operator with a lower bound $C$ there exists its self-adjoint 
extension, which is also semibounded with a lower bound $C'$ for an arbitrary $C'<C$, see~\cite[Satz~43]{jvN29}. In 
addition, as a footnote, von Neumann conjectured that it is even possible to take $C=C'$, i.e., to get a self-adjoint 
extension with the same lower bound. Subsequently, Friedrichs proved the existence of such an extension 
in~\cite[Satz~9]{koF34:MA.a}, and he was even able to specify, under certain specific assumptions on the coefficients, 
the domain of this extension for the second-order Sturm--Liouville differential operator, see~\cite{koF36}. This made 
the extension to be somehow exceptional, and Friedrichs called it as {\it ausgezeichnete Fortzetzung} (an excellent 
extension). His approach was based on an associated quadratic form and ``boundary terms'', which is not, in principle, 
too far from our treatment. The notion of {\it Friedrichs extension} appears probably for the first time in
Freudenthal's work \cite{hF36}, where a limit characterization of this extension was derived for any lower semibounded 
Hermitian linear operator. Actually, this technique turns out to be crucial for many subsequent results (including 
ours). Rellich in~\cite{fR51} provided two alternative characterizations of the Friedrichs extension for the 
second-order Sturm--Liouville differential operator without the conditions imposed on the coefficients by Friedrichs 
instead of which he assumed explicitly that the operator is bounded from below. In particular, he showed that the 
elements of the domain of the Friedrichs extensions behave like a principal solution near the boundary. A similar result 
can also be found in Kalf's paper~\cite{hK78} but this time utilizing Freudenthal's characterization. 

Simultaneously, from the Glazman--Krein--Naimark theorem we know that the domain of any self-adjoint extension of an 
operator can be described by suitable $d$ ``boundary'' conditions, where $d$ is equal to positive and negative 
deficiency indices of the operator, see, e.g., \cite[Theorem~4 in~§18]{maN68}. Zettl and his co-authors showed that 
these are the Dirichlet boundary conditions in the case of the Friedrichs extension of regular ordinary 
differential operators with locally integrable coefficients or in a more general setting, 
see~\cite{hdN.aZ90,qB.gW.aZ22:OaM}. On the other hand, in the singular case the Dirichlet boundary condition at one 
endpoint (or eventually at both endpoints) is not well defined, so another condition is needed in this situation, which 
was treated in~\cite{hdN.aZ92,qB.gW.aZ22:JMAA,qB.gW.aZ23,gmB.jvB09,jvB68,hgK.mkK.aZ86} 
including the description of the Friedrichs extension as a true extension of the minimal operator in~\cite{sY.jS.aZ15}. 
Since, in almost all these cases, a connection between the differential expressions and linear Hamiltonian differential 
systems is used, it is not very surprising that later the Friedrichs extension was solely investigated for operators or 
linear relations associated with these systems itself, see~\cite{mM.aZ00,rSH.pZ10,zZ.qK18,cY.hS21}.

The literature on a discrete analog of this problem is, however, humbler. The Friedrichs extension of the 
Jacobi operator or the second order Sturm--Liouville difference expression was studied 
in~\cite{mB.wdE94,bmB.jsC05,fG.zZ93} and for higher order expressions through the banded symmetric matrices 
in~\cite{oD.pH09}. Furthermore, very recently, Friedrichs extension in the setting of a~linear Hamiltonian difference 
system was investigated in~\cite{sZ.hS.cY23,gR.gX23}. As it is well known that this system can be written as 
a~discrete symplectic system but not vice versa in general, we provide a generalization of the latter results. For 
completeness, we mention that the first attempts of a unification of these results for any even order Sturm--Liouville 
differential and difference expressions through the calculus on time scales were presented 
in~\cite{pZ.pH12,pZ11:IJDSDE}. Our main result (given in Theorem~\ref{T:main}) should not be anyhow surprising as it 
yields the same conclusion as in the continuous case or for Jacobi operators, the latter of which is, in fact, a special 
case of~\eqref{E:Sla}. Nonetheless, it completes this direction by a nontrivial generalization of the above-mentioned 
results, including the linear Hamiltonian difference system. We also note that principal or recessive solutions 
still remain the main tool in the characterization of the Friedrichs extension; although we can find various 
approaches to this characterization in the general theory of linear relations, their application to the specific cases 
mentioned above still seems to be somehow restricted.

The paper is organized as follows. In the next section we introduce the notation used and the basic setting of 
system~\Slaf{\la}{f}, recall a general characterization of the Friedrichs extension in the theory of linear relations 
and the notion of the recessive solution of discrete symplectic systems, and derive several preliminary results. The 
main result is established in Section~\ref{S:main}.

\section{Preliminaries}\label{S:prelim}

Throughout the paper, all matrices are considered over the field of complex numbers $\Cbb$. For $r,s\in\Nbb$ we denote 
by $\Cbb^{r\times s}$ the space of all complex-valued $r\times s$ matrices and $\Cbb^{r\times 1}$ will be abbreviated as 
$\Cbb^r$. In particular, the $r\times r$ \emph{identity} and \emph{zero matrices} are written as $I_r$ and $0_r$, where 
the subscript is omitted whenever it is not misleading (for simplicity, the zero vector is also written as $0$). By 
$e_i$ for $i\in\{1,\dots,n\}$ or $i\in\{1,\dots,2n\}$ we mean the elements of the canonical basis of $\Rbb^n$ or 
$\Rbb^{2n}$, i.e., the columns of $I_n$ or $I_{2n}$. For a given matrix $M\in\Cbb^{r\times s}$
we indicate by $M^*$, $\ker M$, $\rank M$, $M^\dagger$, $M\geq0$, and $M>0$ respectively, its conjugate transpose, 
kernel, rank, the Moore--Penrose generalized inverse, positive semidefiniteness, and positive definiteness. Furthermore, 
we denote by $\Cbb(\oinftyZ)^{r\times s}$ the space of sequences defined on $\oinftyZ$ of complex $r\times s$ matrices, 
where typically $r\in\{n,2n\}$ and $1\leq s\leq2n$. In particular, we write only $\Cbb(\oinftyZ)^r$ in the case $s=1$. 
If $M\in\Cbb(\oinftyZ)^{r\times s}$, then $M(k)\coloneq M_k$ for $k\in\oinftyZ$ and if $M(\la)\in\Cbb(\oinftyZ)^{r\times 
s}$, then $M(\la,k)\coloneq M_k(\la)$ for $k\in\oinftyZ$ with $M_k^*(\la)\coloneq[M_k(\la)]^*$. If 
$M\in\Cbb(\oinftyZ)^{r\times s}$ and $N\in\Cbb(\oinftyZ)^{s\times p}$, then $MN\in\Cbb(\oinftyZ)^{r\times p}$, where 
$(MN)_k\coloneq M_k N_k$ for $k\in\oinftyZ$. We put $\big[z_k \big]_{k=m}^n\coloneq z_n - z_m$. We also adopt a common
notation that $2n$-vector-valued sequences or solutions of \Slaf{\la}{f} are denoted by small letters, typically
$z=\msmatrix{x\\ u}$ with $n$-vector valued components, while $2n\times m$ matrix-valued solutions are denoted by
capital letters, typically $Z=\msmatrix{X\\ U}$ with $n\times m$ matrix-valued components. For completeness, we note
that any solution of \eqref{E:Sla} can be easily seen as a solution of \Slaf{0}{f}.

Finally, the square summability is defined via the semi-inner product
 \begin{equation*}%\label{E:semi-inner.def}
  \innerP{z}{u}\coloneq\sum_{k=0}^\infty z_k^*\,\Ps_k\,u_k \qtextq{and the induced semi-norm}
  \normP{z}\coloneq\sqrt{\innerP{z}{z}}
 \end{equation*}
with respect to the weight matrices $\Ps_k$ specified in~\eqref{E:Sc.Ps.intro}, i.e., we restrict our attention to the
space
 \begin{equation*}%\label{E:ltp.def}
  \ltp=\ltp(\oinftyZ)\coloneq\{z\in\Cbb(\oinftyZ)^{2n}\mid \normP{z}<\infty\}
 \end{equation*}
and, subsequently, to the corresponding Hilbert space
 \begin{equation*}\label{E:tltp.def}
  \tltp=\tltp(\oinftyZ)\coloneq \ltp\big/\big\{z\in\Cbb(\oinftyZ)^{2n}\mid \ \normP{z}=0\big\}
 \end{equation*}
consisting of equivalence classes, which are denoted by $\zclass$.

In the most general setting, a {\it linear relation} $\pzcT$ is defined as a linear subspace of the Cartesian product
of two vector spaces $\pzcX$ and $\pzcY$. We focus only on the case when $\pzcX=\pzcY$ and it is a Hilbert space 
$\pzcH$, which provides suitable tools to study nondensely defined operators via their graphs. For a~deeper insight in 
this theory we refer to~\cite{jB.sH.hsvdS20} and in a~connection with discrete symplectic systems 
to~\cite{slC.pZ15,pZ.slC16,pZ.slC22}. We recall that the {\it domain} of a linear relation 
$\pzcT\subseteq\pzcH\times\pzcH$ is defined as 
 \begin{equation*}
  \dom\pzcT\coloneq\big\{\pzcz\in\pzcH\mid \exists \pzcf\in\pzcH \text{ such that } \{\pzcz,\pzcf\}\in\pzcT\big\}
 \end{equation*}
and the {\it adjoint} relation of $\pzcT$ as  
 \begin{equation*}
  \pzcT^*\coloneq\big\{\{\pzcy,\pzcg\}\in\pzcH^2\mid 
               \inner{\pzcz}{\pzcg}-\inner{\pzcf}{\pzcy}=0\ \text{for all }\ \{\pzcz,\pzcf\}\in \pzcT\big\}.
 \end{equation*}
The linear relation $\pzcT$ is said to be ({\it semi)bounded from below by $c\in\Rbb$}, i.e., $\pzcT\geq c$, if
 \begin{equation}\label{E:semibounded.T}
  \inner{\pzcf}{\pzcz}\geq c\,\inner{\pzcz}{\pzcz} \qtext{for all $\{\pzcz,\pzcf\}\in\pzcT$}
 \end{equation}
and, in particular, {\it nonnegative}, i.e., $\pzcT\geq0$, if
 \begin{equation*}
  \inner{\pzcf}{\pzcz}\geq 0 \qtext{for all $\{\pzcz,\pzcf\}\in\pzcT$.}
 \end{equation*}
The largest $c$ satisfying~\eqref{E:semibounded.T} is said to be the {\it lower bound} of $\pzcT$. In that case the 
linear relation is necessarily symmetric and it has equal defect numbers, which guarantees the existence of its 
self-adjoint extension(s). In particular, in the case of equal deficiency indices there exists the Friedrichs
extension $\pzcT_F$ of $\pzcT$ as defined in \cite[Section~5.3]{jB.sH.hsvdS20}, which can be characterized as follows,
see~\cite[Corollary~5.3.4]{jB.sH.hsvdS20}

\begin{theorem}\label{T:TF.general}
 Let $\pzcT$ be a semibounded linear relation in $\pzcH^2$. Then $\big\{\pzcz,\pzcf\big\}\in\pzcT_F$ if and only if
 $\{\pzcz,\pzcf\}\in\pzcT^*$ and there exists a sequence $\big\{\{\pzcz_n,\pzcf_n\}\big\}_{n=1}^\infty\in\pzcT$
 such that
  \begin{equation*}\label{E:T_F.limit}
   \pzcz_n\to\pzcz \qtextq{and} \inner{\pzcz_n}{\pzcf_n}\to\inner{\pzcz}{\pzcf} \qtext{as $n\to\infty$.}
  \end{equation*}
\end{theorem}

Similarly to the Freudenthal's characterization in the operator case, it can be shown that 
$\{\pzcz,\pzcf\}\in\pzcT_F$ if and only if $\big\{\pzcz,\pzcf\big\}\in\pzcT^*$ and there exist a sequence 
$\big\{\{\pzcz_n,\pzcf_n\}\big\}_{n=1}^\infty\in\pzcT$ such that
 \begin{equation}\label{E:T_F.limit.2}
  \pzcz_n\to\pzcz \qtextq{and} \inner{\pzcz_n-\pzcz_m}{\pzcf_n-\pzcf_m}\to0 \qtext{as $n,m\to\infty$,}
 \end{equation}
see also \cite{sH.aS.hsvdS.hW07:JOT,sH04,sH.mK.hsvdS97:JOP,eaC73:AMS}. 

The following hypothesis summarizes the basic assumptions concerning system~\eqref{E:Sla.noref} or \Slaf{\la}{f} with 
the special linear dependence on the spectral parameter. These systems can be determined either by the pair of
coefficient matrices $\{\Sc,\Vc\}$ or by the pair $\{\Sc,\Ps\}$ and there is no difference between these two approaches
as the matrices $\Ps_k$ and $\Vc_k$ are mutually connected via the equalities $\Ps_k=\Jc\.\Sc_k\.\Jc\.\Vc^*_k\.\Jc$ and
$\Vc_k=-\Jc\,\Ps_k\,\Sc_k$. Furthermore, Hypothesis~\ref{H:basic} yields that system~\eqref{E:Sla} can be written by 
using the matrices $\Sbb_k(\la)\coloneq \Sc_k+\la\.\Vc_k$, which satisfy the symplectic-type equality 
$\Sbb_k^*(\bla)\.\Jc\.\Sbb_k(\la)=\Jc$. This guarantees the existence of a unique solution of any initial value problem 
associated with~\eqref{E:Sla}.

\begin{hypothesis}\label{H:basic}
 A number $n\in\Nbb$ and the pair of matrix-valued sequences $\Sc\in\Cbb(\oinftyZ)^{2n\times2n}$ and 
 $\Wc\in\Cbb(\oinftyZ)^{n\times n}$ are given such that
  \begin{equation*}%\label{E:H.basic}
   \Sc_k^*\Jc\Sc_k=\Jc \qtextq{and} \Wc_k^*=\Wc_k>0 \qtext{for all $k\in\oinftyZ$.}
  \end{equation*}
 The matrices $\Sc_k$ admit the $n\times n$ block decomposition
  \begin{equation*}
   \Sc_k=\mmatrix{\Ac_k & \Bc_k\\ \Cc_k & \Dc_k}
  \end{equation*}
 for all $k\in\oinftyZ$ and the matrix-valued sequences 
 $\Ps,\Vc\in\Cbb(\oinftyZ)^{2n\times2n}$ are defined as
  \begin{equation*}
   \Ps_k\coloneq\mmatrix{\Wc_k & 0\\ 0 & 0}
   \qtextq{and}
   \Vc_k\coloneq-\Jc\,\Ps_k\,\Sc_k \qtext{for all $k\in\oinftyZ$.}
  \end{equation*}
\end{hypothesis}

Our main result combines tools from the spectral theory of linear relations and from the oscillation theory
of~\eqref{E:Sla.noref} with $\la\in\Rbb$, where an important role is played by a ``special'' type of matrix-valued 
solutions, for which we need the following notions.

\begin{definition}
 Let $\nu\in\Rbb$ be fixed and Hypothesis~\ref{H:basic} be satisfied. A $2n\times n$ matrix-valued solution
 $Z(\nu)\in\Cbb(\oinftyZ)^{2n\times n}$ of~\Sla{\nu} is said to be a {\it conjoined basis} if $\rank Z_k(\nu)=n$ and
 $Z_k^*(\nu)\.\Jc\.Z_k(\nu)=0$ for some (and hence for any) $k\in\oinftyZ$. Two conjoined bases
 $Z(\nu),\tZ(\nu)\in\Cbb(\oinftyZ)^{2n\times n}$ of~\Sla{\nu} are said to be {\it normalized} if
 $Z_k^*(\nu)\.\Jc\.\tZ_k(\nu)=I$ for some (and hence for any) $k\in\oinftyZ$.
\end{definition}

A~comprehensive treatise on the qualitative theory of discrete symplectic systems can be found in the recent book
\cite{oD.jE.rSH19}. Our notion of recessive solutions in Definition~\ref{D:recessive} follows the traditional concept 
introduced in~\cite{oD00:EJQTDE} and its generalization was studied in~\cite{pS.rSH15:LAA,pS.rSH17}.

\begin{definition}\label{D:recessive}
 Let $\nu\in\Rbb$ be fixed and Hypothesis~\ref{H:basic} be satisfied. A conjoined basis 
 $\tZ(\nu)=\msmatrix{\tX(\nu)\\ \tU(\nu)}$ of system~\Sla{\nu} is said to be a {\it recessive solution} if, for large
 $k\in\oinftyZ$, the matrix $\tX_k(\nu)$ is nonsingular, it holds $-\tX_k^{-1}(\nu)\.\Bc_k\.\tX_{k+1}^{*-1}(\nu)\geq0$ 
 and simultaneously $\lim_{k\to\infty} X_k^{-1}(\nu)\.\tX_k(\nu)=0$ for any conjoined basis  
 $Z(\nu)$ normalized with $\tZ(\nu)$, i.e., such that $Z_k^*(\nu)\.\Jc\.\tZ_k(\nu)\equiv I$.
\end{definition}

Note that the recessive solution is determined uniquely up to a right multiple by a constant nonsingular $n\times n$
matrix. However, not every system~\eqref{E:Sla.noref} possesses a recessive solution. Its existence can by guaranteed by 
two additional assumptions as we show in Theorem~\ref{T:L1}. Let $\nu\in\Rbb$ be fixed. System~\Sla{\nu} is said to be 
{\it nonoscillatory} if there exists $M\in\oinftyZ$ such that it is {\it disconjugate} on $[M,N+1]_\sZbb$ for
every $N\in[M,\infty)_\sZbb$, i.e., the matrix-valued solution $Z(\nu)\in\Cbb(\oinftyZ)^{2n\times n}$ 
determined by the initial condition $Z_{N+1}(\nu)=\msmatrix{0\\ -I}$ satisfies
 \begin{equation}\label{E:kernel.condition}
  \ker X_{k}(\nu)\subseteq\ker X_{k+1}(\nu) \qtextq{and} -X_{k+1}(\nu)\.X_k^\dagger(\nu)\.\Bc_k\geq0
 \end{equation}
for all $k\in[M,N]_\sZbb$, see~\cite[Theorem~2.41]{oD.jE.rSH19}. In the opposite case, system~\Sla{\nu} is called {\it 
oscillatory}. The nonoscillatory behavior implies that every conjoined basis $Z(\nu)$ of~\Sla{\nu} satisfies 
condition~\eqref{E:kernel.condition} for all $k\in\oinftyZ$ large enough and, consequently, the kernel of $X_k(\nu)$ is 
eventually constant. In addition, we say that system~\Sla{\nu} is {\it disconjugate on $\oinftyZ$} if it is 
nonoscillatory with $M=0$.

System~\eqref{E:Sla.noref} is ({\it completely}) {\it controllable} on a~discrete interval $[N,\infty)_\sZbb$ if for 
any nontrivial finite discrete subinterval $[K,M]_\sZbb\subset[N,\infty)_\sZbb$ the trivial solution $z(\la)\equiv0$ is 
the only solution of \eqref{E:Sla.noref} with $x_k(\la)=0$ for all $k\in[K,M]_\sZbb$, i.e., the subsystem
 \begin{equation}\label{E:control}
  0=\Bc_k\.u_{k+1} \quad\&\quad u_k=\Dc_k\.u_{k+1}, \quad k\in[K,M-1]_\sZbb,
 \end{equation}
has only the trivial solution, see also equations~\eqref{E:Sla.1}--\eqref{E:Sla.2} below. This happens, e.g., when 
$\Bc_k$ is invertible for all $k\in[N,\infty)_\sZbb$. Note that the subsystem in~\eqref{E:control} does not involve 
$\la$, so the controllability can be seen as a global property of system~\eqref{E:Sla.noref} independent of $\la$. 
System \eqref{E:Sla.noref} is {\it eventually controllable} if there exists $N\in\oinftyZ$ such that it is completely 
controllable on $[N,\infty)_\sZbb$. This property together with the eventually constant kernel $X_k(\nu)$ mentioned 
above implies the invertibility of $X_k(\nu)$ for all $k\in\oinftyZ$ large enough, i.e., if system \Sla{\nu} is 
nonoscillatory and eventually controllable, then for every conjoined basis $Z(\nu)$, there exists $N\in\oinftyZ$ such 
that $X_k(\nu)$ is invertible and $-X_{k+1}(\nu)\.X_k^{-1}(\nu)\.\Bc_k\geq0$ for all $k\in[N,\infty)_\sZbb$.

The following result is a time-reversed analogue of~\cite[Theorem~3.1]{oD00:EJQTDE} and
\cite[Theorem~2.66]{oD.jE.rSH19}. We omit its proof because it can be done in the same way as in the mentioned 
references.

\begin{theorem}\label{T:L1}
 Let Hypothesis~\ref{H:basic} hold and $\nu\in\Rbb$ be such that system \Sla{\nu} is nonoscillatory and eventually
 controllable. Then \Sla{\nu} possesses a recessive solution $\tZ=\msmatrix{\tX\\ \tU}\in\Cbb(\oinftyZ)^{2n\times n}$, 
 which can be equivalently characterized by the condition
  \begin{equation*}
   \lim_{k\to\infty} \la_{\min} \Big(-\sum^k_{j=k_0} \tX_j^{-1}(\nu)\.\Bc_j\.\tX_{j+1}^{*-1}(\nu) \Big)=\infty,
  \end{equation*}
 where $\la_{\min}$ stands for the smallest eigenvalue of the matrix indicated and $k_0\in\oinftyZ$ is large enough. 
\end{theorem}

The maximal linear relation is defined as in~\eqref{E:Tmax.def}, while the minimal linear relation displayed
in~\eqref{E:Tmin.explicit} is defined as the closure of the {\it pre-minimal linear relation} $T_0$, which consists of
$\{\zclass,\fclass\}\in\Tmax$ such that $\hz_0=0$ and $\hz_k=0$ for all $k\in\oinftyZ$ large enough and a suitable
representative $\hz\in\zclass$. The following hypothesis guarantees that the minimal linear relation can be written as
in~\eqref{E:Tmin.explicit}. It is called as the {\it strong Atkinson condition} or {\it definiteness condition} and  
it is a classical assumption in the Weyl--Titchmarsh theory for differential or difference equations. In addition, it 
is equivalent to the fact that for any $\{\zclass,\fclass\}\in\Tmax$ there is a~unique $\hz\in\zclass$ such that 
$\mL(\hz)_k=\Ps_k\,f_k$ for all $k\in\oinftyZ$, see \cite[Theorem~5.2]{slC.pZ15}. Since the latter equality 
is also independent of the choice of a representative of $f\in\fclass$, we may write only $\{z,f\}\in\Tmax$ whenever 
the hypothesis is satisfied.

\begin{hypothesis}[Strong Atkinson condition]\label{H:definiteness}
 Hypothesis~\ref{H:basic} is satisfied, and a number $\nu\in\Cbb$ and a finite interval
 $\IzD\coloneq[a,b]_\sZbb\subseteq\oinftyZ$ exist such that every nontrivial solution $z(\nu)\in\Cbb(\oinftyZ)^{2n}$
 of system~\Sla{\nu} satisfies $\sum_{k\in\IzD} z_k^*(\la)\.\Ps_k\.z_k(\la)>0$.
\end{hypothesis}

Note that the strong Atkinson condition is independent of the choice of $\la\in\Cbb$, i.e.,
Hypothesis~\ref{H:definiteness} means that $\sum_{k\in\IzD} z_k^*(\la)\.\Ps_k\.z_k(\la)>0$ is satisfied for all
nontrivial solutions of \eqref{E:Sla.noref} for any $\la\in\Cbb$, see e.g. \cite[Lemma~2.1]{pZ23}. Alternatively, this 
condition is equivalent to the fact that the trivial solution is the only solution of~\eqref{E:Sla.noref} such that 
$\sum_{k\in\IzD} z_k^*(\la)\.\Ps_k\.z_k(\la)=0$. Such systems are also said to be {\it definite} on the discrete 
interval $\oinftyZ$.

In the next part we aim to connect the disconjugacy of~\eqref{E:Sla.noref} on $\oinftyZ$ and the boundedness from below 
of $\Tmin$. Due to the special block structure of $\Sc_k$ and $\Ps_k$ described in Hypothesis~\ref{H:basic}, 
system~\eqref{E:Sla.noref} can be written as the pair of equations
 \begin{align}
  &x_k(\la)=\Ac_k\.x_{k+1}(\la)+\Bc_k\.u_{k+1}(\la) \label{E:Sla.1}\\
  &u_k(\la)=\Cc_k\.x_{k+1}(\la)+\Dc_k\.u_{k+1}(\la)+\la\.\Wc_k\.x_k(\la). \label{E:Sla.2}
 \end{align}
Then a sequence $z\in\Cbb(\oinftyZ)^{2n}$ is said to be {\it admissible} if it satisfies equation~\eqref{E:Sla.1}, 
which does not involve the parameter $\la$ explicitly, i.e., the space of all admissible sequences 
of~\eqref{E:Sla.noref} is independent of $\la$. In the next lemma, we introduce an quadratic functional associated 
with~\eqref{E:Sla} and describe its connection to the inner product $\innerP{\cdot}{\cdot}$, from which we will derive 
the dependence of the disconjugacy on $\la$.
\enlargethispage{0.5mm}

\begin{lemma}\label{L:quadratic.functional}
 Let $\la\in\Cbb$ be arbitrary, Hypothesis~\ref{H:basic} be satisfied, and for any
 $z=\msmatrix{x\\ u}\in\Cbb(\oinftyZ)^{2n}$ define the quadratic functional
  \begin{equation*}
   \Fc_\la(z)\coloneq -\sum_{k\in\oinftyZ} \Big\{x_{k+1}^*\.\Cc_k^*(\la)\.\Ac_k\.x_{k+1}
                         +2\re\!\big(\.x_{k+1}^*\.\Cc_k^*(\la)\.\Bc_k\.u_{k+1}\big)
                           +u_{k+1}^*\.\Dc_k^*(\la)\.\Bc_k\.u_{k+1}\Big\},
  \end{equation*}
 where $\Cc_k(\la)\coloneq \Cc_k+\la\.\Wc_k\.\Ac_k$ and $\Dc_k(\la)\coloneq \Dc_k+\la\.\Wc_k\.\Bc_k$.
 If $z\in\Cbb(\oinftyZ)^{2n}$ is an admissible sequence of~\eqref{E:Sla.noref}, it reduces to
  \begin{align}
   \Fc_\la(z)&=\sum_{k\in\oinftyZ} \big(u_k-\Cc_k(\la)\.x_{k+1}-\Dc_k(\la)\.u_{k+1}\big)^*x_k
                   +\big[u_k^*\.x_k\big]_{k=0}^\infty\notag\\
   &=\sum_{k\in\oinftyZ} \big(u_k-\Cc_k\.x_{k+1}-\Dc_k\.u_{k+1}\big)^*x_k
     +\big[u_k^*\.x_k\big]_{k=0}^\infty-\la\.\innerP{z}{z}\label{E:quadratic.functional.admissible}
  \end{align}
 and, furthermore, if $z$ solves \Slaf{\nu}{f} for some $\nu\in\Cbb$ then
  \begin{equation*}
   \Fc_\la(z)=(\bar{\nu}-\la)\.\innerP{z}{z}+\innerP{z}{f}+\big[u_k^*\.x_k\big]_{k=0}^\infty.
  \end{equation*}
 Especially, when $\{\zclass,\fclass\}\in T_0$, we obtain
  \begin{equation*}
   \Fc_\la(z)=\innerP{z}{f}-\la\.\innerP{z}{z}=\innerP{z}{f-\la\.z}
  \end{equation*}
 for any $z\in\zclass$ and $f\in\fclass$. 
\end{lemma}

It was shown in~\cite[Theorem~3.1]{rSH.pZ12} and~\cite[Theorem~2.41]{oD.jE.rSH19} that the disconjugate property of
system~\Sla{\nu} on $[M,N+1]_\sZbb$ is equivalent with the positivity of the associated quadratic functional
 \begin{equation}\label{E:quadratic.functional.finite}
   -\sum_{k=M}^N \Big\{x_{k+1}^*\.\Cc_k^*(\nu)\.\Ac_k\.x_{k+1}
                         +2\.\re\!\big(x_{k+1}^*\.\Cc_k^*(\nu)\.\Bc_k\.u_{k+1}\big)
                           +u_{k+1}^*\.\Dc_k^*(\nu)\.\Bc_k\.u_{k+1}\Big\}
 \end{equation}
for any admissible $z=\msmatrix{x\\ u}\in\Cbb([M,N+1]_\sZbb)^{2n}$ with $x_M=0=x_{N+1}$ and $x\neq0$. Recall that the 
space of all admissible sequences for~\eqref{E:Sla.noref} is independent on $\la$, so
from~\eqref{E:quadratic.functional.admissible} we get that for any admissible $z\in\Cbb(\oinftyZ)^{2n}$ it holds
 \begin{equation}\label{E:F_nu.F_la}
  \Fc_\la(z)=\Fc_\nu(z)+(\nu-\la)\.\innerP{z}{z}.
 \end{equation}
Therefore, the nonnegativity (or positivity) of $\Fc_\nu(z)$ for some $\nu\in\Rbb$ implies the same property for all 
$\la<\nu$. Subsequently, upon combining with Hypothesis~\ref{H:definiteness}, we get the following corollary, whose
second part is a simple consequence of Theorem~\ref{T:L1}.

\begin{corollary}\label{C:disconjugacy}
 Let Hypothesis~\ref{H:definiteness} be satisfied and $\nu\in\Rbb$ be such that system~\Sla{\nu} is disconjugate on
 $\oinftyZ$ and eventually controllable. Then system~\eqref{E:Sla.noref} is disconjugate on $\oinftyZ$ and possess a 
 recessive solution for any $\la\leq\nu$.
\end{corollary}

The last part of Lemma~\ref{L:quadratic.functional} shows that the boundedness from below of $\Tmin$ is closely 
connected to the nonnegativity of $\Fc_\la(\cdot)$ or, in fact, with the disconjugacy of~\Sla{\nu} on $\oinftyZ$ as 
stated in the next theorem.

\begin{theorem}\label{T:Tmin.bounded}
 Let Hypothesis~\ref{H:definiteness} be satisfied and $\nu\in\Rbb$ be such that system~\Sla{\nu} is disconjugate on
 $\oinftyZ$. Then, $\Tmin$ is bounded from below by a lower bound $c\geq\nu$ or equivalently $\Tmin-\la\.I$ is
 bounded from below by $c-\la>0$ for all $\la<\nu$. Consequently, the deficiency indices of $\Tmin$ satisfy 
 $d_+(\Tmin)=d_-(\Tmin)=d_\mu(\Tmin)$ for all $\mu\in\Cbb\stm[\nu,\infty)$ and any self-adjoint extension of $\Tmin$ is
 bounded from below.
\end{theorem}

\begin{proof}
 Since $\Tmin=\overline{T_0}$, it suffices to show that $T_0$ is bounded from below. From the definition of $T_0$ and
 the positivity of~\eqref{E:quadratic.functional.finite} on any subinterval $[M,N+1]_\sZbb\subset\oinftyZ$ with $M=0$
 we get $\innerP{z}{f-\nu\.z}=\Fc_{\nu}(z)>0$ for all $z\in\dom T_0$. This shows that $T_0$ is bounded from below
 by $c\geq\nu$ or equivalently the boundedness of $\Tmin-\la\.I$ from below by $c-\la>0$ for all $\la<\nu$. The second 
 part of the statement follows immediately from \cite[Proposition~1.4.6]{jB.sH.hsvdS20} and 
 \cite[Proposition~5.5.8]{jB.sH.hsvdS20}.
\end{proof}

The second part of Theorem~\ref{T:Tmin.bounded} shows that the disconjugacy of \Sla{\nu} on $\oinftyZ$ guarantees
the existence of a self-adjoint extension of $\Tmin$, which is possible if and only if the positive and negative 
deficiency indices $d_+(\Tmin)$ and $d_-(\Tmin)$, respectively, are equal, see~\cite[Corollary, p.~34]{eaC73:AMS}. 
This equality $d_+(\Tmin)=d_-(\Tmin)\eqcolon d$ can be alternatively interpreted under Hypothesis~\ref{H:definiteness} 
so that systems~\eqref{E:Sla.noref} and~\Sla{\bla} possess the same number $d$ of linearly independent square summable 
solutions for any $\la\in\Cbb\stm\Rbb$, see~\cite[Corollary~5.12]{slC.pZ15}. If, in addition, there is a number 
$\nu\in\Rbb$ such that system~\Sla{\nu} has the same number $d$ of linearly independent square summable solutions, then 
all self-adjoint extensions of $\Tmin$ admit the following characterization, see \cite[Theorem~3.3 
and~Remark~3.4]{pZ.slC16}. This statement turns out to be yet another crucial ingredients in the proof of our main 
result, see Lemma~\ref{L3}.

\begin{theorem}\label{T:intro.sadj.GKN.Sla}
 Let Hypothesis~\ref{H:definiteness} be satisfied and assume that
  \begin{enumerate}[leftmargin=10mm,topsep=2mm,label={{\normalfont{(\roman*)}}}]
   \item both systems \Sla{i} and \Sla{-i} possess $d$ linearly independent square summable solutions;
   \item there exists $\nu\in\Rbb$ such that also system \Sla{\nu} possess $d$ linearly independent square 
         summa\-ble solutions, which are denoted as (suppressing the argument $\nu$) $\Th^{[1]},\dots,\Th^{[d]}$ and
         arranged so that the $2(d-n)\times2(d-n)$ leading principal submatrix of the $d\times d$ matrix 
         $\Ups\coloneq\Th^*_0\,\Jc\,\Th_0$ has a full rank, where
         $\Th_k\coloneq(\Th_k^{[1]},\dots,\Th_k^{[d]})$ for $k\in\oinftyZ$.
  \end{enumerate}
 Then, a linear relation $T\subseteq \tltp\times\tltp$ is a~self-adjoint extension of the minimal linear relation 
 $\Tmin$ if and only if there exist matrices $M\in\Cbb^{d\times 2n}$ and $L\in\Cbb^{d\times 2(d-n)}$ such that
  \begin{equation}\label{E:T.intro.matrices.M.L.def}
   \rank(M,L)=d,\quad M\Jc M^*-L\,\Ups_{2(d-n)\times 2(d-n)}\,L^*=0,
  \end{equation}
 and 
  \begin{equation}\label{E:T.intro.sa-e}
   T=\Bigg\{\{z,f\}\in\Tmax\mid M z_0-L\msmatrix{(\Th^{[1]},z)_{\infty}\\ \vdots\\ 
            (\Th^{[2(d-n)]},z)_{\infty}}=0\Bigg\},
  \end{equation}
 where $\Ups_{2(d-n)\times 2(d-n)}$ is the $2(d-n)\times 2(d-n)$ principal leading submatrix of $\Ups$ and  
 $(\Th^{[j]},z)_{\infty}\coloneq \lim_{k\to\infty} \Th^{[j]*}_k\.\Jc\.z_k$ for $j=1,\dots,2(d-n)$ exist due to 
 the square summability of both sequences.
\end{theorem}

\section{Main result}\label{S:main}

At this moment we have presented all preliminary results needed to establish our main result. Its proof is based on the 
following three lemmas. In the first lemma, we show that $x_0=0$ for every $\{z,f\}\in T_F$. In the second lemma, we 
prove that the columns of a recessive solution of~\eqref{E:Sla.noref} with $\la<\nu$ belong to the domain of the 
Friedrichs extension of $\Tmin$, which is denoted as $T_F$. Thereafter, in the third lemma, we show that a certain 
linear relation determined by a (part of a) recessive solution is self-adjoint. 

\begin{lemma}\label{L2a}
 Let Hypothesis~\ref{H:definiteness} be satisfied and $\nu\in\Rbb$ be such that system~\Sla{\nu} is disconjugate on
 $\oinftyZ$ and eventually controllable. Then $x_0=0$ for any $\{z,f\}\in T_F$.
\end{lemma}
\begin{proof}
 The assumptions guarantee that the Friedrichs extension of the minimal linear relation $\Tmin$ exists by
 Theorem~\ref{T:Tmin.bounded}. Condition~\eqref{E:T_F.limit.2} together with Hypothesis~\ref{H:basic} means that 
 $\{z,f\}\in T_F$ if and only if $\big\{z,f\big\}\in\Tmax$ and there exists a sequence  
 $\big\{\{z^{[n]},f^{[n]}\}\big\}_{n=1}^\infty\in\Tmin$ such that
  \begin{equation*}
   \lim_{n\to\infty}\normP[\big]{z^{[n]}-z}
    =\lim_{n\to\infty} \sum_{k=0}^\infty (x^{[n]}_k-x_k)^*\.\Wc_k\.(x^{[n]}_k-x_k)=0
  \end{equation*}
 and, by using Lemma~\ref{L:quadratic.functional},
  \begin{equation*}\label{E:dom.T_F.limit.22}
   \innerP{z^{[n]}-z^{[m]}}{f^{[n]}-f^{[m]}}=\Fc_0(z^{[n]}-z^{[m]})
     -\big[\big(u_k^{[n]}-u_k^{[m]}\big)^*\.\big(x_k^{[n]}-x_k^{[m]}\big)\big]_{k=0}^\infty\to0
  \end{equation*}
 as $n,m\to\infty$. Since $\Wc_k>0$ for all $k\in\oinftyZ$ by Hypothesis~\ref{H:basic}, the first condition yields that
 $\lim_{n\to\infty} \Wc_0\.(x^{[n]}_0-x_0)=0$. Since $x^{[n]}_0=0$ for all $n\in\Nbb$ by~\eqref{E:Tmin.explicit}, it
 follows that also $x_0=0$ for every $\{z,f\}\in T_F$.
\end{proof}

\begin{lemma}\label{L2}
 Let Hypothesis~\ref{H:definiteness} be satisfied and $\nu\in\Rbb$ be such that system~\Sla{\nu} is disconjugate on
 $\oinftyZ$ and eventually controllable. Then, for any $\la\leq\nu$, all the columns of the recessive solution 
 $\tZ(\la)$ belong to $\ltp$ and their trivializations at $0$ given as
  \begin{equation*}
    \ttzj(\la)=\begin{cases}
          0,& k\in[0,a]_\sZbb,\\
          \tzj(\la),& k\in[b+1,\infty)_\sZbb
         \end{cases}
  \end{equation*}
 belong to $\dom T_F$ for all $j=1,\dots,n$, where $a,b$ are the endpoints of the discrete interval $\IzD$ from 
 Hypothesis~\ref{H:definiteness}.
\end{lemma}
\begin{proof}
 For better clarity, the proof is divided into three steps, which concerns with the behavior of the columns of 
 $\tZ(\la)$ in a~neighborhood of $\infty$. More precisely, by using a recessive solution of~\eqref{E:Sla.noref} we 
 construct a sequence of $2n\times n$ matrix-valued solutions of this system, which converges to $\tZ(\la)$ (the first 
 step). These solutions give rise to pairs $\big\{\ttzjm,\ttfjm\big\}\in T_0-\la\.I$ (the second step), which 
 satisfy both limit conditions in~\eqref{E:T_F.limit.2} from the characterization of $T_F$ (the third step).

 {\it Step 1.} Since $\la\leq\nu$, the system~\eqref{E:Sla.noref} is disconjugate on $\oinftyZ$ and eventually 
 controllable as well by Corollary~\ref{C:disconjugacy}. Thus, it possesses a recessive solution 
 $\tZ(\la)\in\Cbb(\oinftyZ)^{2n\times n}$ by Theorem~\ref{T:L1} and we denote its columns as
 $\tz^{[1]}(\la),\dots,\tz^{[n]}(\la)$. The disconjugacy of~\eqref{E:Sla.noref} implies, in addition, that 
 $-\tX_k^{-1}(\la)\.\Bc_k\.\tX_{k+1}^{*-1}(\la)\geq0$, and so for the matrix-valued sequence 
 $\La_k\coloneq \sum_{j=0}^{k-1} -\tX_j^{-1}(\la)\.\Bc_j\.\tX_{j+1}^{*-1}(\la)$ we 
 have $\lim_{k\to\infty} \La_k^{-1}=0$ by Theorem~\ref{T:L1}. Let us define the so-called associated dominant solution 
 of~\eqref{E:Sla.noref} as
  \begin{equation*}
   \hX_k(\la)\coloneq \tX_k(\la)\.\La_k \qtextq{and}
   \hU_k(\la)\coloneq \tU_k(\la)\.\La_k+\tX_k^{*-1}.
  \end{equation*}
 Then $\tZ(\la)$ and $\hZ(\la)$ form normalized conjoined bases of~\eqref{E:Sla.noref}, the matrices $\hX_k(\la)$ are 
 eventually nonsingular by the controllability of \eqref{E:Sla.noref}, and $\hX_k^{-1}(\la)\.\tX_k(\la)$ is Hermitian. 
 For a fixed $m\in\oinftyZ$ large enough we define (suppressing the ``dependence'' on $\la$)
  \begin{equation*}
   \Xm_k\coloneq \tX_k-\hX_k\.\hX^{-1}_m\.\tX_m \qtextq{and}
   \Um_k\coloneq \tU_k-\hU_k\.\hX^{-1}_m\.\tX_m, \quad  k\in\oinftyZ.
  \end{equation*}
 Then $\Zm$ is a solution of~\eqref{E:Sla.noref} as a linear combination of two solutions of this system with 
  \begin{equation*}
   \Zm_m=\mmatrix{0\\ \tU_m-\hU_m\.\hX^{-1}_m\.\tX_m}=\mmatrix{0\\ -\hX_m^{*-1}}.
  \end{equation*} 
 Furthermore, it holds
  \begin{align*}
   &\Xm_k=\tX_k-\tX_k\.\La_k\.\La_m^{-1}\.\tX^{-1}_m\.\tX_m=\tX_k\big[I-\La_k\.\La_m^{-1}\big]\to\tX_k,\\
   &\Um_k=\tU_k-\big[\tU_k\.\La_k+\tX_k^{*-1}\big]\.\La_m^{-1}\.\tX_m^{-1}\.\tX_m
         =\tU_k\big[I-\La_k\.\La_m^{-1}\big]-\tX_k^{*-1}\.\La_m^{-1}\to \tU_k
  \end{align*}
 as $m\to\infty$. 
 
 {\it Step 2.} Let
  \begin{equation*}
   \zjm=\mmatrix{\xjm\\ \ujm}\coloneq \Zm\.e_j
  \end{equation*}
 and
  \begin{equation*}
   \tzjm_k=\mmatrix{\txjm_k\\[0.5mm] \tujm_k}\coloneq \begin{cases}
                                          \zjm_k, & k\in[0,m]_\sZbb,\\
                                          0, & k\in[m+1,\infty)_\sZbb,
                                         \end{cases}
   \qtextq{and}
    \tfjm_k\coloneq\begin{cases}
                  0,& k\neq m,\\
                  \msmatrix{-\Wc_m^{-1}\.\hX_m^{*-1}\.e_j\\ 0}, & k=m,
                 \end{cases}
  \end{equation*}
 with $\Wc_m$ being the $n\times n$ left upper block of $\Ps_m$. Then obviously $\tzjm,\tfjm\in\ltp$, the 
 sequence $\tzjm$ is admissible, and by a direct calculation we can verify that 
 $\mL(\tzjm)_k=\Ps_k(\la\.\tzjm_k+\tfjm_k)$ for all $k\in\oinftyZ$, i.e., $\big\{\tzjm,\tfjm\big\}\in\Tmax-\la\.I$ for 
 any $m\in\oinftyZ$. In addition, according to the Patching lemma established in \cite[Lemma~3.1]{pZ.slC16}, we have 
 also $\big\{\ttzjm,\ttfjm\big\}\in\Tmax-\la\.I$ with
  \begin{equation*}
   \ttzjm_k=\begin{cases}
             0,& k\in[0,a]_\sZbb,\\
             \tzjm_k,& k\in[b+1,\infty)_\sZbb
            \end{cases}
   \qtextq{and}
   \ttfjm_k=\begin{cases}
             0,& k\in[0,a]_\sZbb,\\
             \tfjm_k,& k\in[b+1,\infty)_\sZbb,
            \end{cases}
  \end{equation*}
 which yields that $\big\{\ttzjm,\ttfjm\big\}\in T_0-\la\.I$ for all $m>b+1$.

 {\it Step 3.} Now, we show that $\big\{\ttzjm,\ttfjm\}\in T_0-\la\.I\subseteq\Tmax-\la\.I$ is such that, for all 
 $j\in\{1,\dots,n\}$, it satisfies $\ttzjm\to\tzj$ as $m\to\infty$ and simultaneously
  \begin{equation}\label{E:L2.inner.step4}
   \inner{\ttf^{[j,m]}-\ttf^{[j,\ell]}}{\ttz^{[j,m]}-\ttz^{[j,\ell]}}\to0 \qtext{as $m,\ell\to\infty$.}
  \end{equation}
 So, without loss of generality, let $\ell>m>b+1$. Then, by a direct calculation, we get
  \begin{align*}
   \innerP{\ttf^{[j,m]}-\ttf^{[j,\ell]}}{\ttz^{[j,m]}-\ttz^{[j,\ell]}}
    &=\sum_{k=0}^\infty \big(\ttf_k^{[j,m]}-\ttf_k^{[j,\ell]}\big)^*\.\Ps_k\big(\ttz_k^{[j,m]}-\ttz_k^{[j,\ell]}\big)\\
    &\hspace*{-20mm}=\sum_{k=b+1}^\infty 
          \Big(\ttf_k^{[j,m]*}\Ps_k\.\ttz_k^{[j,m]}+\ttf_k^{[j,\ell]*}\Ps_k\.\ttz_k^{[j,\ell]}
            -\ttf_k^{[j,m]*}\.\Ps_k\.\ttz_k^{[j,\ell]}-\ttf_k^{[j,\ell]*}\.\Ps_k\.\ttz_k^{[j,m]}\Big)\\
    &\hspace*{-20mm}=\tf_m^{[j,m]*}\Ps_m\.\tz_m^{[j,m]}+\tf_\ell^{[j,\ell]*}\Ps_\ell\.\tz_\ell^{[j,\ell]}
       -\tf_m^{[j,m]*}\.\Ps_m\.\tz_m^{[j,\ell]}-\tf_\ell^{[j,\ell]*}\.\Ps_\ell\.\tz_\ell^{[j,m]}\\
    &\hspace*{-20mm}=0+0-\mmatrix{-\Wc_m^{-1}\.\hX_m^{*-1}\.e_j\\ 0}^*\mmatrix{\Wc_m & 0\\ 0 & 0}\.
    \mmatrix{\tX_m-\hX_m\.\hX^{-1}_\ell\.\tX_\ell\\ \tU_m-\hU_m\.\hX^{-1}_\ell\.\tX_\ell}-0\\
    &\hspace*{-20mm}=\e_j^*\.\hX_m^{-1}\.\Wc_m^{-1}\.\Wc_m\.\big(\tX_m-\hX_m\.\hX^{-1}_\ell\.\tX_\ell\big)\.e_j\\
    &\hspace*{-20mm}=\e_j^*\.\big(\hX_m^{-1}\tX_m-\hX^{-1}_\ell\.\tX_\ell\big)\.e_j,
  \end{align*}
 where we used the special block structure of $\Ps_k$ and the facts that $\tz_\ell^{[j,m]}=0$ and 
 $\tX_m^{[m]}=0=\tX_\ell^{[\ell]}$. Thus, \eqref{E:L2.inner.step4} holds due to the definition of the recessive 
 solution. Furthermore, from Theorem~\ref{T:Tmin.bounded} we know that $\Tmin-\la\.I$ is bounded from below by a lower 
 bound $c-\la>0$, so
  \begin{equation*}
   \innerP{\ttf^{[j,m]}-\ttf^{[j,\ell]}}{\ttz^{[j,m]}-\ttz^{[j,\ell]}}\geq
    (c-\la)\.\innerP{\ttz^{[j,m]}-\ttz^{[j,\ell]}}{\ttz^{[j,m]}-\ttz^{[j,\ell]}},
  \end{equation*}
 which together with the previous conclusion implies that also
  \begin{equation*}
   \innerP{\ttz^{[j,m]}-\ttz^{[j,\ell]}}{\ttz^{[j,m]}-\ttz^{[j,\ell]}}\to0
  \end{equation*}
 as $m,\ell\to\infty$, i.e., the sequence $\big\{\ttzjm\big\}_{m=0}^\infty$ is a Cauchy sequence in $\ltp$. The
 definiteness condition guarantees that each $\ttzjm$ gives rise to a unique equivalence class, so each column $\tzj$ 
 of the recessive solution belongs to $\ltp$. Consequently, $\ttzj$ belongs to the domain of the Friedrichs 
 extension of $\Tmin-\la\.I$, which is equal to $T_F-\la\.I$ as shown in~\cite[Identity~(5.3.4)]{jB.sH.hsvdS20}.
 Therefore, $\ttzj\in\dom T_F$ for all $j\in{1,\dots,n}$ and the proof is complete.
\end{proof}

In the last lemma, we prove that the linear relation
 \begin{equation}\label{E:U.def}
  U\coloneq\Big\{\{z,f\}\in\Tmax\mid x_0=0 \text{ \ and\ } \lim_{k\to\infty} z_k^*\.\Jc\.\tz_k^{[i_j]}=0
     \text{ for all $j=1,\dots,d-n$}\Big\}
 \end{equation}
is a self-adjoint extension of $\Tmin$, where $\tz^{[i_j]}=\tz^{[i_j]}(\la)$ is the $i_j$-th column of the recessive 
solution $\tZ(\la)$ being such that $\tz^{[i_j]}_m=e_{i_j}$ for a suitable $m\in\oinftyZ$, arbitrary $\la<\nu$, and 
certain indices $i_j$ with $j=1,\dots,d-n$ specified in the proof of Lemma~\ref{L3}, see equality~\eqref{E:L3.ij.def}. 
This is done by a~construction of the matrix $\Ups$ mentioned in Theorem~\ref{T:intro.sadj.GKN.Sla}. Note that we 
do not emphasize the dependence of $U$ on $\la$ in~\eqref{E:U.def}, because we will show in Theorem~\ref{T:main} that
it is only formal in the present context and $U$ represents a Friedrichs extension of $\Tmin$ for all suitable $\la$.

\begin{lemma}\label{L3}
 Let assumptions of~Lemma~\ref{L2} hold and $\la<\nu$ be arbitrary. The linear relation $U$ defined 
 in~\eqref{E:U.def} is a self-adjoint extension of $\Tmin$.
\end{lemma}
\begin{proof}
 Since $\la<\nu$, it follows from Theorem~\ref{T:Tmin.bounded} that system~\eqref{E:Sla.noref} possesses 
 $n\leq d\leq 2n$ linearly independent square summable solutions, so the first assumption of
 Theorem~\ref{T:intro.sadj.GKN.Sla} is satisfied and a~self-adjoint extension of $\Tmin$ exists.
 The proof is therefore completed by, at first, constructing a~matrix-valued solution
 $\Th(\la)\in\Cbb(\oinftyZ)^{2n\times d}$ of system~\eqref{E:Sla.noref} consisting of $d$ linearly independent square
 summable solutions satisfying the second assumption of Theorem~\ref{T:intro.sadj.GKN.Sla}, and thereafter constructing
 matrices $M$ and $L$ satisfying the conditions in~\eqref{E:T.intro.matrices.M.L.def} such that the corresponding
 self-adjoint extension displayed in~\eqref{E:T.intro.sa-e} is equal to $U$.

 From Lemma~\ref{L2} we know that all $n$ columns of a recessive solution $\tZ(\la)$ are square summable and, without
 loss of generality, we may assume that $\tX_m(\la)=I$ for some $m\in\oinftyZ$. We complete these solutions with the
 remaining $d-n$ linearly independent square summable solutions of~\eqref{E:Sla.noref}, which can be taken as the
 columns of some $\hZ(\la)\in\Cbb(\oinftyZ)^{2n\times(d-n)}$. For simplicity, we suppress the dependence on $\la$ in
 the rest of the proof. If we put $\hhZ \coloneq \hZ-\tZ\.\hX_m\in\Cbb(\oinftyZ)^{2n\times(d-n)}$, then it
 solves~\eqref{E:Sla.noref} and $\hhZ_m=\msmatrix{0\\ \hhhU_m}$. Since $\hhX_m=0$, it follows that $\rank \hhU_m=d-n$
 and one can easily deduce that, after an appropriate constant multiple of $\hhZ$, there are indices
 $i_1,\dots,i_{d-n}$ such that
  \begin{equation}\label{E:L3.ij.def}
   \mmatrix{e_{i_1}^*\\ \vdots\\ e_{i_{d-n}}^*}\.\hhU_m=I_{d-n}.
  \end{equation}

 Then we can build the $2n\times d$ matrix-valued solution $\Th=\big(\Th^{[1]}\ \ \Th^{[2]}\ \ \Th^{[3]}\big)$ mentioned
 above from $\Th^{[2]}\coloneq \hhZ$ and the blocks
  \begin{gather*}
   \Th^{[1]}\coloneq \tZ\mmatrix{e_{i_1}^*\\ \vdots\\ e_{i_{d-n}}^*}^{\!*}\in\Cbb(\oinftyZ)^{2n\times(d-n)}
   \qtextq{and}
   \Th^{[3]}\coloneq \tZ\mmatrix{e_{s_1}^*\\ \vdots\\ e_{s_{2n-d}}^*}^{\!*}\in\Cbb(\oinftyZ)^{2n\times(2n-d)},
  \end{gather*}
 where the indices $i_1,\dots,i_{d-n}\in\{1,\dots,n\}$ correspond to the choice of rows of $\hhU_m$ described 
 in~\eqref{E:L3.ij.def} and $s_{1},\dots,s_{2n-d}\in\{1,\dots,n\}\stm\big\{i_1,\dots,i_{d-n}\big\}$ are all the 
 remaining indices. To justify this choice we need to show that the $2(d-n)\times2(d-n)$ leading principal submatrix of
 the $d\times d$ matrix $\Ups\coloneq\Th^*_0\,\Jc\,\Th_0$ has a full rank. Since $\la\in\Rbb$, it holds
 $\Ups=\Th^*_m\,\Jc\,\Th_m$ by the Wronskian-type identity, see \cite[Identity~(3.4)]{pZ.slC22}, and the submatrix can
 be decomposed as
  \begin{equation*}
   \Ups_{2(d-n)\times2(d-n)}
    =\mmatrix{\Ups_{2(d-n)\times2(d-n)}^{[1,1]} & \Ups_{2(d-n)\times2(d-n)}^{[1,2]}\\[2mm]
             \Ups_{2(d-n)\times2(d-n)}^{[2,1]} & \Ups_{2(d-n)\times2(d-n)}^{[2,2]}}
    =\mmatrix{\Th_m^{[1]*}\.\Jc\Th_m^{[1]} & \Th_m^{[1]*}\.\Jc\Th_m^{[2]}\\[2mm]
              \Th_m^{[2]*}\.\Jc\Th_m^{[1]} & \Th_m^{[2]*}\.\Jc\Th_m^{[2]}}.
  \end{equation*}
 Then $\Ups_{2(d-n)\times2(d-n)}^{[1,1]}=0$ as it is a submatrix of $\tZ^*_m\.\Jc\.\tZ_m$, which is zero by the
 definition of the recessive solution (it has to be a conjoined  basis), while $\tX_m=I$ and $\hhX_m=0$ yield
 $\tZ_m^*\.\Jc\.\hhZ_m=\hhU_m$, so
  \begin{equation*}
   \Ups_{2(d-n)\times2(d-n)}^{[1,2]}=-\Ups_{2(d-n)\times2(d-n)}^{[2,1]}
    =\msmatrix{e_{i_1}^*\\ \vdots\\ e_{i_{d-n}}^*}\.\hhU_m=I_{d-n}
  \end{equation*}
 by~\eqref{E:L3.ij.def} and~\eqref{E:matrix.J.def}. In addition, $\hhX_m=0$ implies also
 $\Ups_{2(d-n)\times2(d-n)}^{[2,2]}=\hhZ_m^*\.\Jc\.\hhZ_m=0$. Therefore, 
  \begin{equation}\label{E:Ups.submatrix.I}
   \rank\Ups_{2(d-n)\times2(d-n)}=\rank\mmatrix{0 & I_{d-n}\\ -I_{d-n} & 0}=2(d-n),
  \end{equation}
 i.e., the matrix-valued solution $\Th$ satisfies the assumption (ii) of Theorem~\ref{T:intro.sadj.GKN.Sla}.

 It remains to express the linear relation $U$ from~\eqref{E:U.def} as in~\eqref{E:T.intro.sa-e}. If we put
  \begin{equation*}
   M\coloneq \mmatrix{I_n & 0\\ 0 & 0} \qtextq{and} 
   L\coloneq \mmatrix{0 & 0\\ I_{d-n} & 0},
  \end{equation*}
 then we can verify by a simple calculation that these matrices satisfy the conditions 
 in~\eqref{E:T.intro.matrices.M.L.def}. Simultaneously, the equality
  \begin{equation*}
   0=M\.z_0-L\msmatrix{(\Th\.e_{1},z)_{\infty}\\ \vdots\\ (\Th\.e_{2(d-n)},z)_{\infty}}
    =\mmatrix{x_0\\ 0_{d-n}}-\mmatrix{0_n\\ (\Th\.e_1,z)_{\infty}\\ \vdots\\ (\Th\.e_{d-n},z)_{\infty}}
    =\mmatrix{x_0\\ (\Th\.e_{1},z)_{\infty}\\ \vdots\\ (\Th\.e_{d-n},z)_{\infty}}
  \end{equation*}
 utilized in~\eqref{E:T.intro.sa-e} is equivalent to the pair of conditions $x_0=0$ and 
 $\lim_{k\to\infty} z_k^*\.\Jc\.\tz_k^{[i_j]}=0$, because  $\Th\.e_{j}=\tz^{[i_j]}$ for all $j\in\{1,\dots,d-n\}$. 
 Therefore, the linear relation $U$ is a~self-adjoint extension of $\Tmin$.
\end{proof}

Now, upon combining the preceding lemmas and the self-adjointness of the Friedrichs extension of $\Tmin$ and of the 
linear relation $U$, we obtain the main result showing that $T_F=U$. 

\begin{theorem}\label{T:main}
 Let Hypothesis~\ref{H:definiteness} be satisfied and $\nu\in\Rbb$ be such that system~\Sla{\nu} is disconjugate on
 $\oinftyZ$ and eventually controllable. Then, for any $\la<\nu$, the linear relation $U$ defined in~\eqref{E:U.def} 
 is the Friedrichs extension of $\Tmin$, i.e., 
  \begin{equation*}
   T_F=\Big\{\{z,f\}\in\Tmax\mid x_0=0 \text{ \ and\ } \lim_{k\to\infty} z_k^*\.\Jc\.\tz_k^{[i_j]}(\la)=0
     \text{ for all $j=1,\dots,d-n$}\Big\}.
  \end{equation*}
 In particular, if system \eqref{E:Sla.noref} is in the limit point case (i.e., $d=n$), then 
  \begin{equation*}
   T_F=\Big\{\{z,f\}\in\Tmax\mid x_0=0\Big\},
  \end{equation*}
 while in the limit circle case (i.e., $d=2n$), we have
  \begin{equation*}
   T_F\coloneq\Big\{\{z,f\}\in\Tmax\mid x_0=0 \text{ \ and\ } \lim_{k\to\infty} z_k^*\.\Jc\.\tz_k^{[j]}(\la)=0
     \text{ for all $j=1,\dots,n$}\Big\}.
  \end{equation*}
\end{theorem}
\begin{proof}
 We recall that the given assumptions guarantee the existence of the Friedrichs extension of $\Tmin$ by 
 Theorem~\ref{T:Tmin.bounded}. We already know that the linear relation $U$ is a self-adjoint extension of $\Tmin$ by 
 Lemma~\ref{L3}. Now, let $\big\{z,f\big\}\in T_F$ be arbitrary. Then also $\big\{z,f\}\in T_F^*$ and $x_0=0$ by 
 Lemma~\ref{L2a}. Let $\ttz^{[i_1]},\dots,\ttz^{[i_{d-n}]}$ be as in Lemma~\ref{L2}, in particular, 
 $\ttz^{[i_1]},\dots,\ttz^{[i_{d-n}]}\in \dom T_F$ and they coincide with $\tz^{[i_1]},\dots,\tz^{[i_{d-n}]}$ for all
 $k$ large enough. The self-adjointness of $T_F$ yields
  \begin{equation*}
   0=\innerP{f^{[j]}}{z}-\innerP{\ttz^{[i_j]}}{f}
    =\big[\big(\ttz^{[i_j]},z\big)_k\big]_{k=0}^\infty=\lim_{k\to\infty} z_k^*\.\Jc\.\ttz_k^{[i_j]}
    \qtext{for all $j\in\{1,\dots,d-n\}$,}
  \end{equation*}
 where $f^{[1]},\dots,f^{[d-n]}$ are such that 
 $\big\{\ttz^{[i_1]},f^{[1]}\big\},\dots,\big\{\ttz^{[i_{d-n}]},f^{[d-n]}\big\}\in T_F$. This shows that 
 $\big\{z,f\big\}\in U$, which means $T_F\subseteq U$. However, the self-adjointness of $U$ implies also 
 the opposite inclusion $U=U^*\subseteq T_F^*=T_F$ and, therefore, $T_F=U$. The rest of the statement is a simple 
 consequence of this result in the case when $d=n$ and $d=2n$, respectively.
\end{proof}

\begin{remark}\label{R:overdetermined}
 \begin{enumerate}[leftmargin=10mm,topsep=2mm,label={{\normalfont{(\roman*)}}}]
  \item\label{R:overdetermined.i}
   Besides the more general setting in our treatise, which could be easily applied also in the continuous case, 
   Theorem~\ref{T:main} is a discrete analogue of \cite[Theorem~4.2]{zZ.qK18}, \cite[Theorem~3.1]{rSH.pZ10}, and 
   \cite[Theorem~13]{mM.aZ00} concerning the Friedrichs extension of {\it operators} associated with the linear 
   Hamiltonian differential system. It is worth noting that, in contrast to the mentioned results, we neither require 
   any particular number of linearly independent square summable solutions (e.g., the limit circle case as in 
   \cite[Theorem~3.1]{rSH.pZ10}) nor formally overdetermine $T_F$ as in~\cite[Theorem~13]{mM.aZ00}, because we 
   explicitly deal only with $d$ ``boundary''  conditions in accordance with Theorem~\ref{T:intro.sadj.GKN.Sla} --- 
   $n$ of them come from $x_0=0$ by the first part of the proof of Lemma~\ref{L2}, while the remaining $d-n$ specify 
   the behavior at infinity and they are derived from~\eqref{E:L3.ij.def} so that \eqref{E:Ups.submatrix.I} holds, 
   compare with~\cite[Remark after Theorem~12]{mM.aZ00}. 
  
  \item The square summability of the columns of a recessive solution $\tZ(\la)$ established in Lemma~\ref{L2} will be 
   an object of our further research because it seems to be an essential property rooted in the definition of 
   this solution similarly as in the continuous case, see \cite{rSH.pZ18}.  
 \end{enumerate}
\end{remark}

%%%%%%%%%%%%%%%%%%%%%%%%%%%%%%%%%%%%%%%%%%%% SECTION %%%%%%%%%%%%%%%%%%%%%%%%%%%%%%%%%%%%%%%%%%%%%%%%%%%%%%%%%%%%%%%%%

\section*{Acknowledgments}
\addcontentsline{toc}{section}{Acknowledgments}
The author is very grateful to the referee for his/her detailed reading of the manuscript and for various comments and 
suggestions, which improved the presentation of the paper. The author also thanks Roman \v{S}imon Hilscher for
bringing this topic to his attention and for discussions at the beginning of this research. The research was supported
by the Czech Science Foundation
under Grant GA23-05242S.

%%%%%%%%%%%%%%%%%%%%%%%%%%%%%%%%%%%%%%%%%%%% SECTION %%%%%%%%%%%%%%%%%%%%%%%%%%%%%%%%%%%%%%%%%%%%%%%%%%%%%%%%%%%%%%%%%

% \bibliographystyle{my_bibliography_paper_style}
% \bibliography{bibliotheca}

\end{document}